\documentclass[10pt,reqno]{amsart}
\usepackage{amsmath,amssymb,graphicx}

\numberwithin{equation}{section}
\numberwithin{figure}{section}
\numberwithin{table}{section}
\newtheorem{theorem}{Theorem}[section]

\newtheorem{lem}{Lemma}[section]

\newcommand{\ignore}[1]{}{}

\usepackage{amssymb,color}
\definecolor{c20}{rgb}{0.,0.7,0.}
\definecolor{c30}{rgb}{0.,0.,1.}
\definecolor{c40}{rgb}{1,0.1,0.7}
\definecolor{c50}{rgb}{1,0,0}

\def\cK#1{\textcolor{c20}{#1}}
\def\cK#1{#1}

\def\IME#1{\textcolor{c20}{#1}}
\def\IME#1{#1}

\def\jE#1{\textcolor{black}{#1}}

\def\cE#1{\textcolor{c20}{#1}}
\def\aE#1{\textcolor{c20}{#1}}

\def\dE#1{\textcolor{c20}{#1}}
\def\eE#1{\textcolor{c20}{#1}}
\def\eE#1{#1}
\def\hE#1{\textcolor{c20}{#1}}
\def\hE#1{#1}
\def\pE#1{\textcolor{c30}{#1}}
\def\qE#1{\textcolor{c40}{#1}}

\def\pE#1{#1}
\def\qE#1{#1}
\def\cE#1{#1}
\def\dE#1{#1}

\def\aE#1{\textcolor{c20}{#1}}
\def\aE#1{#1}
\def\lE#1{\textcolor{c30}{#1}}
\def\lE#1{#1}
\def\kE#1{\textcolor{c30}{#1}}
\def\kE#1{#1}
\def\JE#1{\textcolor{c30}{#1}}
\def\JE#1{#1}

\def\jj#1{\textcolor{c30}{#1}}

\def\jj#1{#1}

\def\enk#1{{#1}}

\def\Me#1{\textcolor{c40}{#1}}
\def\Me#1{{#1}}

\def\LA{\aE{\widetilde{\lambda}_n}}

\newcommand{\BQN}{\begin{eqnarray}}
\newcommand{\EQN}{\end{eqnarray}}
\newcommand{\BQNY}{\begin{eqnarray*}}
\newcommand{\EQNY}{\end{eqnarray*}}

\def\Pr{\mathbb{P}}
\newcommand{\pk}[1]{\Pr \left(#1\right) }

\newcommand{\R}{\!I\!\!\mathbb{R}}

\newcommand{\N}{\mathbb{N}}
\newcommand{\inr}{\in \mathbb{R}}

\newcommand{\prooftheo}[1]{ \textsc{Proof of Theorem} \ref{#1} }

\newcommand{\prooflem}[1]{\textsc{Proof of Lemma} \ref{#1}}

\newcommand{\QED}{\hfill $\Box$}
\newcommand{\IF}{\infty}
\newcommand{\BL}{\begin{lem}}
\newcommand{\EL}{\end{lem}}
\newcommand{\nelem}[1]{{Lemma \ref{#1}}}

\newcommand{\netheo}[1]{{Theorem \ref{#1}}}

\newcommand{\todis}{\stackrel{d}{\to}}
\newcommand{\toprob}{ \stackrel{p}{\to}}
\newcommand{\ldot}{,\ldots,}
\newcommand{\E}[1]{\mathbb{E}\{#1\}}

\def\limit{\lim_{t\to \IF}}

\def\CM{C^*}
\def\CMF{\widetilde{C}}
\newcommand{\abs}[1]{\lvert #1 \rvert}

\begin{document}

\begin{center}

\vspace*{1cm}

\MakeUppercase{\bf Tail asymptotics  of randomly weighted 
 large risks }

\large \normalsize

\bigskip

\textsc{Alexandru V. Asimit}\footnote{Corresponding author. Phone: +44(0)2070405282. Fax: +44(0)2070408572.}

\textit{Cass Business School, City University, London EC1Y 8TZ, \\
United Kingdom.} \texttt{E-mail}: \texttt{asimit@city.ac.uk}

\bigskip

\textsc{Enkelejd Hashorva\footnote{This work is partially supported by the Swiss National Science Foundation Grant 200021-134785.}}

\textit{Department of Actuarial Science, Faculty of Business and Economics,\\
University of Lausanne, UNIL-Dorigny 1015 Lausanne, Switzerland.}\\ \texttt{E-mail}: \texttt{Enkelejd.Hashorva@unil.ch}

\bigskip

\textsc{Dominik Kortschak\footnote{This work was partially supported  by the
the MIRACCLE-GICC project and the Chaire d'excellence "Generali - Actuariat responsable: gestion des risques naturels et changements climatiques."}}

\textit{Universit\'e de Lyon, F-69622, Lyon, France; Universit\'e Lyon 1, Laboratoire SAF, EA 2429, Institut de Science Financi\`ere et d'Assurances, 50 Avenue Tony Garnier, F-69007 Lyon, France}\\ \texttt{E-mail}:  \texttt{Kortschakdominik@gmail.com}

\bigskip

\today{} 

\bigskip

\end{center}

\hspace{1in}

\begin{quote}
\noindent
\textbf{Abstract.}
\pE{Tail asymptotic probabilities for linear combinations of randomly weighted order statistics are approximated under various assumptions. One key assumption is the asymptotic independence for all risks, and thus, it is not surprising that the maxima represents the most influential factor when one investigates the tail behaviour of our considered risk aggregation, which for example, can be found in the reinsurance market. This extreme behaviour confirms the ``one big jump" property that has been vastly discussed in the existing literature in various forms whenever the asymptotic independence is present. An illustration of our results together with a specific application are explored under the assumption that the underlying risks follow the multivariate Log-normal distribution.}

\noindent
\textit{Keywords and phrases}:
Davis-Resnick tail property; Extreme value distribution; Max-domain of attraction; Mitra-Resnick model; 
Risk aggregation.
\end{quote}

\section{Introduction}
Consider positive dependent random variables (or risks) $X_i,i=1, \ldots, n$ and let $X_{1,n}\geq\ldots\geq X_{n,n}$ be the corresponding upper order statistics. We investigate the asymptotic tail behaviour probability for linear combinations of order statistics, $L(\textbf{C}) = \sum_{i=1}^k C_i X_{i,n},\;1\leq k\leq n,$ where  the $C_i$'s are random deflators/weights and  $C_1>0$.

Studying the tail probability for such order statistics has multiple financial and insurance applications (for example, see Hashorva, 2007, Ladoucette and Teugels, 2006, Asimit and Jones, 2008 a and b, Jiang and Tang, 2008 and Li and Hashorva, 2013). \pE{All of these papers have studied the extreme behaviour of the $L(\textbf{C})$, where $\textbf{C}$ was assumed to be a constant random vector.} 
The tail asymptotics of $L(\textbf{1})$ has been investigated in many recent contributions such as Asmussen and Rojas-Nandayapa (2008), Chen and Yuen (2009), Mitra  and  Resnick  (2009), Foss  and Richards  (2010), Asmussen~\textit{et al.} (2011),   Kortschak (2012),  Hashorva (2013), Embrechts~\textit{et al.} (2014), Hashorva~\textit{et al.} (2014). 
\pE{It has been seen that the assumption of constant $C_i$'s represents a popular setting considered in the recent past, where the asymptotic tail probabilities of some linear combinations of order statistics have been obtained. Obviously, randomising the $C_i$'s is a more challenging problem to be studied, which is the main purpose of this paper. It is also of interest to recognise situations in which the randomisation represents a problem of interest. This is the case if one is interested in a more accurate risk aggregation, where the time value of the money is introduced in the model. That is, this popular risk evaluation takes into account not only the amount of claim, but also the time \qE{when} the claim occurs, and therefore $L(\textbf{C})$ becomes the discounted value of the aggregate risk. Another application that will be detailed in Section~\ref{appl sec}, is given when the $C_i$'s quantify the random proportions paid by the risk holder in the case of its default in payment.}

\pE{We now explain the mathematical framework that will be further assumed in this paper. Under some mild conditions, it is obtained that the most significant contribution to the tail probability of $L(\textbf{C})$ is given by the largest component, i.e., $C_1 X_{1,n}$. This can be explained by the fact that under asymptotic independence, the ``one big jump" property is always present. In other words, as has been observed in the existing literature, the largest value is the most influential factor in risk aggregation. Our proofs are extremely sensitive to the tail behaviour of the individual risks. Therefore, a characterization of the tail distribution of a random variable is necessary, which is a classical result of the Extreme Value Theory. A distribution function (df) $F$ is said to belong to the \textit{Maximum Domain of Attraction} (MDA)\ of a non-degenerate df  $G$, written as $F\in \mathrm{MDA}(G)$,
if there are some $a_{n}>0$ and $b_{n}\inr$ for $n\in \N$ such that \hE{for any constant $x$} $\lim_{n\rightarrow \infty }F^{n}(a_{n}x+b_{n})=G(x),$ where
$G$ is of one of the following three df's:
\begin{eqnarray*}
\begin{array}{lllll}
\text{Fr\'{e}chet:} &  & \Phi _{\alpha }(x)=\exp (-x^{-\alpha })%
 , &  & x>0,\alpha >0; \\
\text{Gumbel:} &  & \Lambda (x)=\exp (-\exp(-x))%
 , &  & -\infty <x<\infty;  \\
\text{Weibull:} &  & \Psi _{\alpha }(x)=\exp (-\abs{x}^{\alpha })%
 , &  & x\leq 0,\alpha >0.
\end{array}
\end{eqnarray*}
We focus on distributions with unbounded support,  i.e., from $\mathrm{MDA}(\Phi_{\alpha })$ and \Me{$\mathrm{MDA}(\Lambda)$}, and therefore only the Fr\'{e}chet and Gumbel cases will be considered. The following section presents our main result. In Section~\ref{appl sec}, we illustrate our findings with an application, while all the proofs are relegated to the last section.}

\section{Main Results}\label{Main Res Sec}

We consider first the Fr\'echet MDA and further include the case where the index $\alpha$ is allowed to be 0, i.e., $X_1$ may exhibit a slowly regularly varying tail. \hE{The mathematical formulation of the tail condition imposed on $X_1$ is given by}
\BQN\label{fresh}
\lim_{t\to \IF} \frac{ \pk{X_1> tx}}{ \pk{X_1> t}}= x^{-\alpha}, \quad \alpha \ge 0.
\EQN
\cE{For more details on regular variations of random variables and vectors see e.g., Jessen and Mikosch (2006).} \\
Throughout the remainder of the paper $\lambda_1 \ldot \lambda_n$ are non-negative constants and  $\widetilde{\lambda}_n:= \sum_{i=1}^n \lambda_i.$

\Me{Some standard notation are used, as well as further explanations, in order to provide a precise meaning of our statements. For two positive functions
$a(\cdot )$ and $b(\cdot )$, we write $a(\cdot )\sim cb(\cdot )$ to mean \qE{asymptotic} equivalence, i.e., $\lim a(\cdot )/b(\cdot )=c$ for some positive
constant $c$. We also denote $\limsup a(\cdot )/b(\cdot )\leq 1$ by $a(\cdot )\lesssim b(\cdot )$.}

\begin{theorem} \label{Th1}
Let $X_1,\ldots,X_n$ be some positive random variables satisfying
\begin{equation}\label{suff1 tail equiv cond}
\lim_{t\to\infty}\frac{\Pr(X_i>t)}{\Pr(X_1>t)}=\lambda_i\in[0,\infty),\;\mbox{for all}\;i\in\{\cE{1},\ldots,n\}.
\end{equation}
Let $\textbf{C}=(C_1 \ldot C_n)$ be a random vector such that $C_1>0$
 is independent of the maximum $X_{1,n}$.
\dE{Suppose that} \eqref{fresh} holds with $\alpha> 0$ and
$\E{C_1^\beta}< \IF$ for some $\beta\in (\alpha,\IF)$. It is assumed that
\begin{equation}\label{suff1 joint tail cond}
\lim_{t\to\infty}\cE{\max_{ 1 \le i < j \le n}} \frac{\Pr(\CMF X_i>t,\CMF X_j>t)}{\Pr(X_1>t)}=0,
\end{equation}
with $\CMF= \max_{2 \le i \le n} |C_i|$. If further there exists a positive constant $\tau$ such that
\begin{eqnarray}\label{Assumption A1}
\Pr\big(\widetilde C X_i > t,  \widetilde C X_j> t, \widetilde C > \tau \big) \ge
\kappa_{ij}\Pr\big(X_i > t/r_{ij},   X_j> t/r_{ij} \big)
\end{eqnarray}
holds for all large $t$ and any two indices $i<j$ in $\{1 \ldot k\}$ with $r_{ij}$ a positive constant, then
\begin{equation}\label{Main Res1}
\Pr\big(L(\textbf{C})>t\big)\sim \lE{\Pr(C_1 X_{1,n}>t)}\sim \Pr(X_1>t) \E{C_1^{\alpha}} \widetilde{\lambda}_n , \quad t\rightarrow\infty.
\end{equation}
\end{theorem}
\textbf{Remark 1:} a) Relation (\ref{Assumption A1}) is clearly satisfied if  $\widetilde C $ is independent of $X_i,i\le k$. Another case
\IME{for which (\ref{Assumption A1}) still holds is}
$\widetilde C $ has a distribution function with positive lower endpoint $\widetilde{\alpha}$ \IME{(take}
$\tau= r_{ij}=\widetilde{\alpha}$). \\
b) \pE{In numerous applications, risks can be of different nature in terms of their tail behaviour, where both light and heavy-tailed risks can be part of the aggregation process. A classical result in the case of independent risks with non-random weights states that the heaviest tail represents the dominant factor in explaining the extreme events of the aggregate risk. This is also the case for our considered model. Indeed, if $X_k$ has heavier tail than $X_1$, i.e.,
\eqref{suff1 tail equiv cond} holds with $\lambda_k=0$, then the result in \eqref{Main Res1} shows that there is no impact of $X_k$ when performing asymptotic evaluations. If $\E{\CMF^\beta} < \IF$ for some $\beta \in (\alpha, \IF)$, then by Lemma \ref{L0} for any $i\not=k$
$$ \frac{\Pr(\CMF X_i>t,\CMF X_k>t)}{\Pr(X_1>t)}\le \frac{\Pr(\CMF X_k>t)}{\Pr(X_1>t)} \to 0, \quad t\to \IF$$
and thus condition \eqref{suff1 joint tail cond}  can be relaxed as follows:
\BQN\label{eqNew}
\lim_{t\to\infty}\cE{\max_{ (i,j)\in E_+} \frac{\Pr(\CMF X_i>t,\CMF X_j>t)}{\Pr(X_1>t)}=0,
\EQN
where $E_+:=\{(i,j): 1 \le i \le j, \lambda_i>0, \lambda_j>0\}$. Note in passing that due to the dependence among the  $C_i$'s, 
\eqref{Assumption A1} is still needed even  if $\lambda_k=0$ for some $k\le n$.
}}\\
\IME{c) }
\jE{It is worth mentioning that (\ref{Assumption A1}) implies} 
\ignore{
\begin{eqnarray*}
\frac{\Pr\big(\widetilde C X_i > t,  \widetilde C X_j> t \big)}{\pk{X_1> t}} & \ge &
\frac{\Pr\big(\widetilde C X_i > t,  \widetilde C X_j> t, \widetilde C > \tau \big)}{\pk{X_1> t}}\\
& \ge &\kappa_{ij}\frac{\Pr\big(X_i > t/r_{ij},   X_j> t/r_{ij} \big)}{\Pr\big(X_1> t/r_{ij}\big)} \frac{\Pr\big(X_1> t/r_{ij}\big)}{\pk{X_1> t}} \to 0 , \quad  t\to \IF.
\end{eqnarray*}
\IME{Consequently, we have}
}
\begin{equation}\label{Assumption 1 res}
\pk{X_i > t,   X_j> t }=o\big( \Pr(X_1> t)\big), \quad  t\to \IF.
\end{equation}

\hE{
Next, we discuss the case in which the random scaling are independent of the portfolio risks, and has the advantage of being able to characterise the tail behaviour \jE{of $L(\textbf{C})$} 
in the presence of slowly variation property of the individual risks.}

\begin{theorem} \label{th3}
Let $X_1,\ldots,X_n$ be some positive random variables satisfying
\pE{\eqref{eqNew}}. Let $\textbf{C}=(C_1 \ldot C_n)$ be a random vector independent of $X_i,i\le n$ with $C_1>0$.
\dE{Suppose that} \eqref{fresh} holds with $\alpha\ge 0$ and
$\E{C_1^\beta}< \IF$ for some $\beta\in (\alpha,\IF)$. If \eqref{suff1 joint tail cond} holds with $\CMF=1$
and further $\eE{\E{\max_{\pE{2} \le i\le n} \enk{\abs{C_i}}^{\gamma}}}\in \pE{(}0, \IF)$
for some $\gamma\in (\alpha,\IF)$, then \eqref{Main Res1} holds and \pE{moreover}
\begin{eqnarray}\label{cle}
\Pr(L(\textbf{C})>t)\sim \sum_{i=1}^n \Pr(C_i X_{i,n}> t), \quad t\rightarrow\infty.
\end{eqnarray}
\end{theorem}

\hE{The Fr\'echet scenario requires a Pareto-like extreme behaviour for the individual risks, and sometimes leads to an overestimate of the extreme events magnitude. Therefore, the Gumbel tail assumption represents a valid alternative, which includes from moderately heavy-tailed distributions, such as Log-Normal, to light-tailed distributions, such as Exponential. We further investigate this scenario for which some background is now provided.}

It is well-known (see Embrechts  et al.\ 1997) that if $F\in \mathrm{MDA}(\Lambda )$, then there exists a positive, measurable function $a(\cdot )$ such that $\bar{F}:=1- F$ satisfies
\begin{equation}\label{GPD Gumbel}
\lim_{t\to \IF}\frac{\bar{F}\left( t+a(t)x\right) }{\bar{F}(t)} = \exp(-x)
\end{equation}
for any $x\inr$. In addition, the latter holds locally uniformly in $x$ (see Resnick, 1987). Recall that the auxiliary function $a(\cdot )$
satisfies $a(t)=o(t)$ and is such that the relation
\begin{equation}
\lim_{t\rightarrow \infty }\frac{a\left( t+a(t)x\right) }{a(t)}=1 \label{a func prop}
\end{equation}
holds locally uniformly in $x$.

If $X_1$ has df $F_1$ satisfying \eqref{GPD Gumbel}, then the asymptotic tail of $X_1$ can be determined under some weak conditions on the random scaling factor $C_1$.
 In the \jE{sequel} we shall \hE{only consider} the case \hE{in which} $C_1$ is bounded,  i.e., its df  has some finite upper endpoint $\omega \in (0,\infty)$. \hE{Specifically, the following two settings are investigated in this paper:}
\begin{enumerate}
\item[i)] \textbf{\underline{Model A}}: Assume that \hE{$\Pr(C_1= \omega)=p\in (0,1]$ and $\Pr(C_1\le \eta)=1-p$ hold for some $\eta \in (0,\omega)$.}
\item[ii)] \textbf{\underline{Model B}}: \hE{For any $x>0$ and some $\gamma \in [0,\IF)$, we have}
\begin{eqnarray}\label{modelB}
\lim_{t\to \IF} \frac{ \pk{C_1> \omega-  x/t}}{\pk{C_1> \omega -  1/t}}= x^\gamma.
\end{eqnarray}
\end{enumerate}
\hE{The case in which $C_1$ has an unbounded upper endpoint,  i.e., $\omega=\infty$, is more complex and less tractable}. If both $C_1$ and
$X_{1,n}$ have Weibullian tails the exact tail asymptotic behaviour of $C_1X_{1,n}$ is derived in Arendarczyk and  D\c{e}bicki (2011).\\
 The paper of Mitra and Resnick (2009) derives the asymptotic tail behaviour of the sum of dependent random variables with Gumbel tails, and their sufficient conditions provide the appropriate framework to elaborate our next result, stated as Theorem~\ref{Th2}. It is worth mentioning that our proof also provides a simplified argumentation of their main result for portfolios consisting of three or more risks.
\JE{In the following we assume without loss of generality that $C_1$ has upper endpoint equal to 1.}
\begin{theorem} \label{Th2}
Let $X_1,\ldots,X_n$ be some positive random variables and \cE{suppose that $X_1$ has df} in the $MDA(\Lambda)$ with infinite right \cE{endpoint} and an auxiliary function $a(\cdot )$ as defined in (\ref{GPD Gumbel}). Let further $\mathbf{C}=(C_1 \ldot C_n)$ be a given random vector such that $C_1>0$
is independent of $X_{1,n}$ and the assumption of Model A or Model B holds \JE{with $\omega=1$}. 
If $(\ref{suff1 tail equiv cond})$ is satisfied and further
\begin{eqnarray}\label{suff2-cond1}
\lim_{t\rightarrow \infty } \max_{1\leq i\neq j\leq n}\frac{\Pr \left(  \CM X_{i}>t, \CM X_{j}>\hE{a}(t)x\right) }{\Pr(C_1 X_1>t)}=0,\quad
 \mathrm{  for\;all}\;x>0
\end{eqnarray}
holds for $\CM:= \max_{1 \le i \le n}\pE{|C_i|}$, then as $t\to \infty$
\begin{eqnarray*}\label{Main Res2}
\Pr(L(\textbf{C})>t)\sim \Pr(C_1 X_{1,n}>t)\sim \dE{\LA \pk{C_1 X_1> t}},
\end{eqnarray*}
provided that
\begin{eqnarray}\label{suff2-cond2}
\lim_{t\rightarrow \infty }\frac{\Pr \left(\CM  X_{i}>L_{ij}a(t),\CM  X_{j}>L_{ij} a(t)\right) }{\Pr(C_1 X_1>t)}=0,\; \mathrm{for\;some}\;L_{ij}>0\;\mathrm{and \;all}\;1\leq i< j\leq n.
\end{eqnarray}
\end{theorem}
\noindent\textbf{\jE{Remark 2}:} In the above theorem, if $C_1$ has upper endpoint $\omega \in (0,\IF)$, then $C_1X_{1,n}$ is in the Gumbel
 MDA with auxiliary function $a_1(t):=\omega a\big(t/\omega\big)$. Clearly, when $\omega=1$, then $a_1(t)=a(t)$. \\
\noindent\textbf{\jE{Remark 3:}} A simpler case is when $C_1 \ldot C_n$ are independent of $X_1 \ldot X_n$, and this setting is addressed \pE{below} in greater details under the framework of Log-Normal risks $X_1 \ldot X_n$.

\section{Applications}\label{appl sec}
\newcommand{\equaldis}{\stackrel{d}{=}}
\pE{In this section, we discuss an application of one of our main results, namely an illustration of the asymptotic approximations found in Theorem~\ref{Th2}. It is further assumed that the individual risks are multivariate Log-Normal distributed. That is, $X_i=\exp\big(\sigma_iZ_i+ \mu_i\big), i\le n,$  where $(Z_1,\ldots,\qE{Z}_n)$ is a multivariate Gaussian distributed random vector with Pearson correlation coefficients $\rho_{ij} \in (-1,1), 1 \le i < j \le n$ and $\sigma_i>0, \mu_i \in \R$ are some given constants. Further, suppose that $X_i,i\le n$ are independent of $C_1$, which is a random variable with $[0,1]$ support. Furthermore, $C_1$ is Beta distributed with positive parameters $\alpha,\beta$. It is well-known that $X_i$'s have distribution functions in the Gumbel \cK{MDA} with scaling functions  $a_i(t)=\sigma_i^2 t/\big(\log(t)-\mu_i\big), i\le n$ (see for example, Embrechts~\textit{et al.}, 1997). By using Lemma~\ref{Th2.0}, it is not difficult to find that $C_1 X_{1,n}$ \qE{has also  df} in the Gumbel MDA.}

\pE{The next step is to show that Theorem~\ref{Th2} is applicable in the current setting, and therefore we only need to establish that condition \eqref{suff2-cond2} holds. In view of our assumptions, we have the stochastic representation
\[(X_i,X_j)\equaldis \left(\exp\left(\sigma_i^2 \sin(\Theta) R+\mu_i\right),\exp\left(\sigma_j^2 \left(\rho_{ij} \sin(\Theta)+\sqrt{1-\rho_{ij}^2}\cos(\Theta) \right)R+\mu_j\right)\right), \quad 1 \le i \not =j < n,\]
where the random angle $\Theta$ is uniformly distributed on $(0,2\pi)$ being independent of the random radius $R>0$. \qE{Note that} $\pk{R>r}=e^{-r^2/2}, r>0$. Assume without loss of generality that $\mu_1>\max_{2 \le i \le n} \mu_i$ and   $\sigma_1\ge \max_{2 \le i \le n} \sigma_i$. 
\qE{Requiring further that $C_i$'s are bounded, say by 1, we have}
$$\Pr \left(\CM  X_{i}>L_{ij}a(t),\CM  X_{j}>L_{ij} a(t)\right) \approx \exp{\left\{-\frac{(\log t)^2}{2\sigma_1^2 \eta_{ij}^2}\right\}}$$
where $\approx$ stands for a logarithmic asymptotic equivalence, and
$$\eta_{ij}=\max_{0\le\theta\le2\pi} \left( \min\left(\sin(\theta),\rho_{ij}\sin(\theta) +\sqrt{1-\rho_{ij}^2}\IME{\cos(\theta)}\right)\right)<1.$$
Note that $f(t)\approx g(t)$ means that $\log f(t)\sim \log g(t)$ is true. Now, \eqref{suff2-cond2} holds since
$\Pr\left(\jj{C_1} X_1>t\right)\approx \exp{\left\{-\frac{(\log t)^2}{2\sigma_1^2}\right\}}$
is satisfied for all large $t$. Recall that $\lambda_i=1 $ if $\sigma_i=\sigma_1, \mu_i=\mu_1$ and otherwise $\lambda_i=0$. Consequently, Theorem~\ref{Th2} holds in our setting.}

\pE{Let us now \qE{discuss} the parametric model considered above. Log-Normal assumption is widely accepted by practitioners as a flexible distribution for modelling individual risks. Choosing an appropriate dependence to model the association among various risks is a difficult tasks, and therefore, the Gaussian dependence is an acceptable choice due to many convenient features, such as the availability of relatively simple estimation and simulation methods. Thus, assuming that a portfolio of $n$ risks is multivariate Log-Normal distributed, represents a reasonable and practical approach. Consider the situation in which the holder of this portfolio, named \textit{insurer}, prefers to transfer the first $k$ largest \qE{claim} amounts, with $1\le k<n$, to a different insurance player, namely \textit{reinsurer}. This risk transfer contract is also known as the \textit{Large Claims Reinsurance (LCR)} (see for example, Ladoucette and Teugels, 2006) and the reinsurer is liable to pay $\displaystyle\sum_{i=1}^k X_{i,n}$, which \qE{might not always be} paid in full due the possibility of default in payment. Therefore, the insurer expects to pay an additional amount (as a result of the default event) of $\displaystyle\sum_{i=1}^k C_i X_{i,n}$, where $0\le C_i=1-RecR_i\le1$ are some random weights with $RecR_i$ being the so-called \textit{recovery rate} corresponding to the $i^{th}$ largest claim. Our assumptions require that $RecR_1$ is Beta distributed with parameters $\beta$ and $\alpha$. Therefore, Theorem~\ref{Th2} tells us that $\Pr\left(\displaystyle\sum_{i=1}^k C_i X_{i,n}>t\right)\sim \Pr\big(C_1X_1>t\big)$, and consequently, the insurer may easily understand the severity of the extreme events associated with the reinsurer default. Specifically, the asymptotic result can be used in approximating tail risk measures such as \textit{Value-at-Risk (VaR)} and \textit{Expected Shortfall (ES)}:
\[ES_p \left(\sum_{i=1}^k C_i X_{i,n}\right)\sim VaR_p \left(\sum_{i=1}^k C_i X_{i,n}\right) \sim VaR_p\big(C_1 X_1\big)\;\;\mbox{as}\;p\uparrow 1,\]
since $C_1 X_{1,n}$ is in the Gumbel MDA (for more details, see Asimit and Badescu, 2010).
Recall that for a generic random variable $Z$, $VaR_p(Z)$ represents the $p^{th}$ quantile and $ES_p(Z):=\E{Z \lvert Z>VaR_p(Z)}$. It can be easily seen that evaluating the extreme events associated with $\displaystyle\sum_{i=1}^k C_i X_{i,n}$ has been drastically reduced via our findings from Theorem~\ref{Th2}.}

\def\kal#1{{\mathcal{ #1}}}

\section{Further Results and Proofs}

\cE{We display  next some lemmas 
which are of some independent interests and then proceed with the proofs of the main results.}

\BL \label{Th2.0}
Let $C_1$ and $X_1$ be two independent positive random variables. Suppose that $C_1$ has upper endpoint $\omega \in (0,\IF)$ and $X_1$
 has df in $MDA(\Lambda)$ with scaling function $a(\cdot)$.
\begin{enumerate}
\item[i)] If $C_1$ obeys Model A, then  \cK{for a function $R(t)$ with $R(t)=o \left(\left(\frac{a(t)}{t}\right)^\xi\right)$, for every $\xi >0$ }
\begin{eqnarray*}
 \pk{C_1 X_1> t \omega} \sim p \pk{X_1> t}\left(1+\cK{R(t)}  \right).
\end{eqnarray*}
\item[ii)]  If  $C_1$ satisfies the assumption \eqref{modelB} of Model B, then
\begin{eqnarray*}
 \pk{C_1 X_1> t \omega} \sim \Gamma(\gamma+1)\Pr\left(C_1 > \omega -  \frac{\omega a(t)}{ t}\right)\pk{X_1> t},
\end{eqnarray*}
 where $\Gamma(\cdot)$ is the Euler gamma function.
\end{enumerate}
\EL

\prooflem{Th2.0} i) The crucial asymptotic result for establishing the claim is the so-called Davis-Resnick tail property of distributions in the \hE{$MDA(\Lambda)$. Namely, by} Proposition~1.1 in Davis and Resnick (1988)  
\begin{eqnarray*}\label{Davis}
\limit  \left(\frac{a(t)}{t}\right)^{\IME{\mu}} \frac{\pk{X_1> K t}}{\pk{X_1> t}}&=& 0
\end{eqnarray*}
holds for any
$\IME{\mu} \inr , K>1$. The latter and the fact that $\eta\in (0,\omega)$ implies the proof. \QED
\ignore{
yield that $\pk{C_1 X_1 > t\omega, C_1< \aE{\eta}} \le \pk{ X_1 > t \omega/\aE{\eta}}=o\left(\left(\frac{a(t)}{t}\right)^\xi \pk{ X_1 > t}\right)$. Thus,
\BQNY
 \pk{C_1 X_1 > t \omega}&=& \pk{C_1 X_1 > t\omega, C_1 =\omega}+ \pk{C_1 X_1 > t\omega, C_1 \le \aE{\eta}}\\
 &\sim & p \pk{ X_1 > t}+o\left(\left(\frac{a(t)}{t}\right)^\xi  \Me{\pk{X_1 > t}}\right) \quad \text{ as } t\to \IF,
 \EQNY
which concludes the first claim.}
\\
ii) Since $C_1^*:=C_1/\omega$ has df with upper endpoint 1 and is regularly varying at 1 with index $\gamma$, the claim follows immediately from Theorem 3.1 in 
Hashorva   et al.\ (2010), and thus the proof is now complete. \QED

{
\BL\label{L0}  Let $C,X$ \hE{and} $Y$ be three  positive random variables such that $C$ is independent from $X,Y$ and $X$ satisfies \eqref{fresh} with some constant
$\alpha \ge 0$. Suppose that \hE{there exists a positive function $h(\cdot)$ such that $\lim_{t\to \IF} h(t)=\lim_{t\to \IF} t/h(t)= \IF$ and}
$\lim_{t\to \IF} \frac{\pk{C> h(t)}}{\pk{X> t}}=0$. If further
\BQN\label{tail}
\lim_{t\to \IF} \frac{\pk{Y> t}}{\pk{X> t}}=0,
\EQN
then we have
\BQN\label{AA}
\lim_{t\to \IF} \frac{\pk{CY> t}}{\pk{X> t}}=0.
\EQN
Furthermore, if $\E{C^{\beta}}< \IF$ for some $\beta> \alpha$, then \eqref{AA} is valid.
\EL
\def\ve{\varepsilon}
\prooflem{L0}
By the assumptions
\BQNY
\frac{\pk{CY> t}}{\pk{X> t}}&=& \frac{\pk{CY> t, C\le   h(t)}+ \pk{CY> t, C >   h(t)}}{\pk{X> t}}\\
&\le& \frac{\pk{CY> t, C\le  h(t)}+ \pk{C >h(t)}}{\pk{X> t}}
\EQNY
holds for all large $t$. \hE{\pE{Denote by} $F$ the df of $C$ and set $h^*(t)=t/h(t)$. \lE{Since $\lim_{t\to \IF} h^*(t)=\lim_{t\to \IF} h(t)=\IF$,
then for any large $M$  we can find $n(M)$ so that for all  $t > n(M)$ we have
$ h^*(t)> M$ and $h(t)> M\enk{.}
$
\cK{Further, by \eqref{tail} for any $\ve >0$ and for some $M'$ (take for simplicity $M'=M$), we have}
$\frac{\pk{Y> t}}{\pk{X> t}} \le \ve, \forall t> M.$ Consequently, for any $c\in (0,1)$ we have
$ h^*(t)/c> M/c> M$ implying
$\pk{Y> h^*(t)/c}/\pk{X> h^*(t)/c} \le \ve, \forall t> n(M).$
}}
The independence of $C$ and $Y$ together with equation \eqref{tail} yield $\Big(\mbox{set}\; G(c):=F\big(h(t) c\big)\Big)$
\BQNY
\frac{\pk{CY> t, C \le   h(t)}}{\pk{CX> t}}&=& \int_{0}^{h(t)}\frac{\pk{Y> t/c}}{\pk{CX> t}}\, d F(c)\\
&=& \frac {1}{\pk{C X> t}}
\int_0^{1}\frac{\pk{Y> h^*(t)/c}}{\pk{X> h^*(t)/c}}\pk{X> h^*(t)/c}\, d \lE{G(c)}\\
&\le & \frac {\varepsilon}{\pk{C X> t}}
\int_0^{1}\pk{X> h^*(t)/c}\, d \lE{G(c)}\\
&= & \varepsilon\frac {\Pr\big(CX> t, C \le h(t)\big)}{\pk{C X> t}} \le \lE{\ve }
\EQNY
\hE{for any  \Me{$t>n(M)$},} and thus \eqref{AA} follows. Next, if $\E{C^{\beta}}< \IF$ for some $\beta> \alpha$, then
the random variable $X$ has heavier tail \Me{than} $C$,  i.e., $\E{X^{\alpha+\varepsilon}}=\infty$ and $\E{(C^{1+ \varepsilon'})^{\alpha+ \varepsilon''}}< \IF$ for some $\varepsilon,\varepsilon', \varepsilon''$ positive with $\varepsilon''> \varepsilon$ implying thus
\begin{eqnarray*}\label{eqCC}
\lim_{t\to \IF} \frac{\pk{C>  \IME{h(t)}}}{\pk{X> t}}=0
\end{eqnarray*}
\IME{since further by the assumption $\lim_{t\to \IF} h(t)/t=0$, hence the proof is complete.}  \QED
%
%

\bigskip
}

In the following we shall need the concept of vague convergence. Let $\{\mu _{n},n\geq 1\}$ be a sequence of measures on a
locally compact Hausdorff space $\mathbb{B}$ with countable base. Then $\mu _{n}\;$converges vaguely to some measure $\mu $, written\ as
$\mu _{n}\overset{v}{\rightarrow }\mu $, if for all continuous functions $f$ with compact support we have
\begin{equation*}
\lim_{n\rightarrow \infty }\int_{\mathbb{B}}f\;\mathrm{d}\mu _{n}=\int_{\mathbb{B}}f\;\mathrm{d}\mu .
\end{equation*}
A thorough background on vague convergence is given in Resnick (1987).

\BL \label{Lm1} \pE{If} the assumptions of \netheo{Th1} are satisfied \pE{with} $\alpha>0$ and $C_i\ge 0$ almost surely for all \hE{$2\leq i\leq n$}, \pE{ then} the following vague convergence
\begin{equation}\label{Th1-step2- vague conv order stats}
\frac{\Pr\big(( C_1 X_{1,n}/t\ldot C_n X_{n,n}/t)\in \cdot\big)}{\Pr( C_1 X_{1,n}>t)}\stackrel{v}{\rightarrow}\mu(\cdot), \quad t \to \IF
\end{equation}
holds on $[0,\infty]\times[M,\infty] \times \cdots \times [M,\IF] \setminus\{(0,0\ldot 0)\}$ for any $M<0$ where the limit measure $\mu$ is given by
\begin{equation*} \label{mu measure}
\mu(dx_1,dx_2\ldot d x_n):=\alpha x_1^{-\alpha-1} dx_1\epsilon_{0}(dx_2) \cdots \epsilon_{0}(dx_n),
\end{equation*}
where $\epsilon_{0}(\cdot)$ denotes the Dirac measure.
\EL

\prooflem{Lm1} First note that by Bonferroni's inequality  for any real $t$ we have
\BQN \label{bon}
\sum_{i=1}^n \Pr(X_i>t)- \sum_{1\leq i<j\leq n} \Pr(X_i>t,X_j>t)\leq \Pr(X_{1,n}>t)\leq \sum_{i=1}^n \Pr(X_i>t).
\EQN
Clearly, equation (\ref{Assumption 1 res}) suggests that $\sum_{1\leq i<j\leq n} \Pr(X_i>t,X_j>t)=o\big(\pk{X_1> t}\big), \quad t\to \IF.$
Hence, \hE{equation} (\ref{suff1 tail equiv cond}) implies that $\lim_{t\to \IF}  \frac{ \Pr( X_{1,n}>t)}{\Pr(X_1>t)}=  \widetilde{\lambda}_n,$ which in turn by Breiman's Lemma (see Breiman, 1965) \enk{yields}
\begin{equation}\label{Th1-step1}
\lim_{t\to \IF}  \frac{ \Pr( C_1 X_{1,n}>t)}{\Pr(X_1>t)}= \E{C_1^\alpha}  \widetilde{\lambda}_n.
\end{equation}
Next, we show the vague convergence only for the first two \qE{largest} order statistics, since  the high dimensional case follows easily by using further the fact that $X_{1,n}\ge X_{2,n} \ge \cdots \ge X_{n,n}$ almost surely.
The above-mentioned convergence of measures holds if the convergence is valid over the following relative compact sets:
\begin{enumerate}
\item[i)] $(x,\infty]\times(y,\infty]$, where $x>0$, $y\geq M$ and $y\neq 0$;
\item[ii)] $\cK{[}0,\infty]\times(y,\infty]$, where $y>0$.
\end{enumerate}
Part i) is now investigated
for which $\Pr(C_1 X_{1,n}>tx,C_2 X_{2,n}>ty)=\Pr(C_1X_{1,n}>tx)$ holds for all $y< 0$ due to the positivity assumption of the $X_i$'s. Consequently, equation~(\ref{Th1-step1}) yields
\[\lim_{t\to\infty}\frac{\Pr\big((C_1 X_{1,n}/t,C_2 X_{2,n}/t)\in (x,\infty]\times(y,\infty]\big)}{\Pr(C_1 X_{1,n}>t)}=x^{-\alpha}=\mu\left((x,\infty]\times(y,\infty]\right).\]
\hE{For any $y>0$, the following is true}
\begin{eqnarray}\label{min}
\Pr(X_{2,n} >y)\leq \sum_{1\leq i<j\leq n} \Pr(X_i>y,X_j>y).
\end{eqnarray}
\hE{Now,}
\begin{eqnarray}\label{Th1-step2-part i}
\lefteqn{\frac{\Pr\big((C_1 X_{1,n}/t,C_2 X_{2,n}/t)\in (x,\infty]\times(y,\infty]\big)}{\Pr(C_1 X_{1,n}>t)}}\nonumber\\
&\leq& \frac{\Pr(C_2 X_{2,n}>t y)}{ \cE{\pk{C_1 X_{1,n}> t}}}\nonumber\\
&\leq& \frac{\Pr(X_1>t y)}{\Pr(X_1>t)} \frac{\Pr(X_1>t)}{\Pr(C_1 X_{1,n}>t)}\sum_{1\leq i<j\leq n} \frac{\Pr(C_2 X_i>t y,C_2 X_j>t y)}{\Pr(X_1>t y)}\nonumber\\
&\sim&  y^{-\alpha}\left(\E{C_1^{\alpha}} \widetilde{\lambda}_n \right)^{-1} \sum_{1\leq i<j\leq n} \frac{\Pr(C_2 X_i>t y,C_2 X_j>t y)}{\Pr(X_1>t y)}\nonumber\\
&\rightarrow& 0=\mu\left((x,\infty]\times(y,\infty]\right) , \quad t\rightarrow\infty,
\end{eqnarray}
\hE{where} the third implication is due to \hE{equations} (\ref{fresh}) and (\ref{Th1-step1}),  
while the fourth implication is a consequence of (\ref{suff1 joint tail cond}). Thus, part i) is fully justified. Finally, part ii) can be
shown in the same manner as displayed in (\ref{Th1-step2-part i}).
\QED

\BL \label{Lm2}
Let us assume that the assumptions of \netheo{Th2} are satisfied such that $C_i\ge 0$ \hE{almost surely for all $2\leq i\leq n$}. Then as $t\to \IF$
\begin{equation}\label{Lemma2- vague conv order stats}
\frac{\Pr\Big(\big( (C_1X_{1,n}- t)/a(t),C_2X_{2,n}/a(t)\ldot C_n X_{n,n}/a(t)\big)\in \cdot\Big)}{\Pr(C_1 X_{1,n}>t)}\stackrel{v}{\rightarrow}\nu(\cdot)
\end{equation}
holds on $[M,\IF] \times [-\IF,\infty]\times \cdots \times [-\IF,\infty]$ for any $M<0$  \hE{with limiting} measure $\nu$ given by
\begin{equation*} \label{mu measure}
\nu(dx_1,dx_2\ldot d x_n):=\exp(- x_1) dx_1\epsilon_{0}(dx_2) \cdots \epsilon_{0}(dx_n).
\end{equation*}
\EL

\prooflem{Lm2}
The proof is similar to that of \nelem{Lm1} and it is sufficient to verify the convergence only over the following compact sets:
\begin{enumerate}
\item[i)] $(x_1,\infty]\times (x_2,\infty]\times \ldots \times (x_n,\infty]$, where $x_1>M$ and $x_i<0$ for all $i\geq 2$;
\item[ii)] $(x_1,\infty]\times (x_2,\infty]\times \ldots \times (x_n,\infty]$, where $x_1>M$ and $x_i>0$ for some $i\geq 2$.
\end{enumerate}
Any set from part i) suggests that 
\[\Pr \big(C_1 X_{1,n}>t+a(t)x_1,C_i X_{i,n}>a(t)x_i, \;\mbox{for all}\; i\geq 2\big)=\Pr\big(C_1 X_{1,n}>t+a(t)x_1\big).\]
\hE{\jE{By Lemma~\ref{Th2.0}, \jE{under the assumptions of Model~A or Model~B the random variable} $C_1 \kE{X_1}$ is in the $MDA(\Lambda)$ with auxiliary function $a(\cdot)$}. \JE{Note that when $\omega \not=1$, then the auxiliary function is not $a(\cdot)$ but $\omega a(t/\omega)$.} Moreover, $C_1 X_{1,n}$ is also in the $MDA(\Lambda)$ with auxiliary function $a(\cdot)$ under the assumptions of Model~A or Model~B,  \kE{ provided that $X_{1,n}$ is in $MDA(\Lambda)$
with auxiliary function $a\jE{(\cdot)}$, which we show next.}}
\kE{We shall show that
\begin{eqnarray}\label{LS}
\lim_{t\to \IF} \frac{\Pr(C_1 X_{1,n}> t)}{\Pr(C_1 X_1> t)}= \widetilde{\lambda}_n
\end{eqnarray}
holds, and thus both $X_{1,n}$ and $C_1 X_{1,n}$ have df in the Gumbel MDA. The fact that $C_1>0$, equations (\ref{suff1 tail equiv cond}) and (\ref{bon}), and the main result of Lemma  \ref{Th2.0} suggest that \eqref{LS} is satisfied as long as
\begin{eqnarray}\label{last eq}
\Pr(C_1 X_i>t,C_1X_j>t)=o\big(\Pr(C_1 X_1>t)\big) , \quad t\to\infty,\;\mbox{for all}\;1\leq i<j \leq n.
\end{eqnarray}
Recall that $a(t)=o(t)$ as $t\to\infty$. The latter and equation (\ref{suff2-cond1}) yield that for large $t$ we have
\begin{eqnarray*}
\Pr(C_1X_i>t,C_1X_j>t)&\leq&\Pr(C^*X_i>t,C^*X_j>t)\\
&\leq&\Pr\big(C^*X_i>t,C^*X_j>a(t)\big)=o\big(\Pr(C_1 X_1>t)\big) , \quad t\to\infty,
\end{eqnarray*}
which justifies (\ref{last eq}).
}
\pE{Consequently,}
\begin{eqnarray*}
\lefteqn{\lim_{t\to \IF}\frac{\Pr\Bigg(\Big(\frac{C_1 X_{1,n}-t}{a(t)},\frac{C_2X_{2,n}}{a(t)},\ldots, \frac{C_n X_{n,n}}{a(t)}\Big)\in (x_1,\infty]\times(x_2,\infty]\times \ldots \times (x_n,\infty]\Bigg)}{\Pr(C_1X_{1,n}>t)}}\\
&=&\lim_{t\to \IF}\frac{\Pr\big(C_1X_{1,n}>t+a(t)x_1\big)}{\Pr(C_1X_{1,n}>t)}\\
&=&e^{-x_1}=\nu\big((x_1,\infty]\times(x_2,\infty]\times \ldots \times (x_n,\infty]\big).
\end{eqnarray*}
For the second part, without loss of generality $x_2>0$ is further assumed. Now,
\begin{eqnarray*}
\lefteqn{\frac{\Pr\Bigg(\Big(\frac{C_1 X_{1,n}-t}{a(t)},\frac{C_2X_{2,n}}{a(t)},\ldots, \frac{C_n X_{n,n}}{a(t)}\Big)\in (x_1,\infty]\times(x_2,\infty]\times \ldots \times (x_n,\infty]\Bigg)}{\Pr(C_1X_{1,n}>t)}}\\
&\leq& \frac{\Pr\Bigg(\Big(\frac{C_1 X_{1,n}-t}{a(t)},\frac{C_2 X_{2,n}}{a(t)}\Big)\in (x_1,\infty]\times(x_2,\infty]\Bigg)}{\pk{C_1X_{1,n}>t}}\\
&\leq& \frac{\pk{C_1X_1>t}}{\pk{C_1X_{1,n}>t}}\frac{\Pr\big( C_1 X_{1,n} > t+ a(t) x_1, C_2 X_{2,n}> a(t) x_2\big)}{\pk{C_1X_1>t}}\\
&\leq& \frac{\Pr\big( \CM X_{1,n} > t+ a(t) x_1,  \CM X_{2,n}> a(t) x_2\big)}{\pk{C_1X_1>t}}\\
&\leq& \sum_{1\leq i\neq j\leq n}\frac{\Pr\big(\CM  X_i>t+a(t)x_1, \CM  X_j>a(t)x_2\big)}{\pk{C_1X_1>t}}\\
&\rightarrow& 0=\nu\big((x_1,\infty]\times(x_2,\infty]\times \ldots \times (x_n,\infty]\big),\quad t\rightarrow\infty
\end{eqnarray*}
follows from (\ref{suff2-cond1}), \cE{\eqref{min}} and the fact that
\begin{eqnarray*}\label{Mmax}
\max_{1\leq i\neq j \leq n} \Pr\big( \CM X_i>t+a(t)x_1, \CM X_j>a(t)x_2\big)=o\big(\Pr(C_1 X_1>t)\big).
\end{eqnarray*}
The latter is justified in few steps. Equation (\ref{a func prop}) yields that $(1-\varepsilon )a(t)\leq a\left( t+a(t)x_1\right) \leq (1+\varepsilon )a(t)$
for any arbitrarily fixed $0<\varepsilon <1$ and all large $t$.
\hE{Recall that $C_1 X_1$ is in the $MDA(\Lambda)$ with scaling function $a(\cdot)$. Consequently,}
\begin{eqnarray*}
\frac{\Pr \left(  \CM X_i>t+a(t)x_1, \CM  X_j>a(t)x_2\right) }{\Pr(C_1 X_1>t)}&\leq &\frac{\Pr \left(  \CM X_i>t+a(t)x_1,  \CM X_j>a\left( t+a(t)x_1\right) \frac{x_2}{1+\varepsilon}\right) }{\Pr(C_1X_1>t)} \\
&\sim &\frac{\Pr \left(  \CM  X_i>t+a(t)x_1,  \CM  X_j>a\left( t+a(t)x_1\right) \frac{x_2}{1+\varepsilon}\right) }{\Pr(C_{1}X_1>t+a(t)x_1)}\exp(-x_1)\\
&\rightarrow& 0, \quad t\rightarrow\infty
\end{eqnarray*}
\hE{holds for any $x_1\inr$ and $x_2>0$}, which is a consequence of relations (\ref{GPD Gumbel}) and (\ref{suff2-cond1}), 
and thus the claim follows. \QED

\prooftheo{Th1}
In the first instance, we assume that $C_i\geq 0$ for all $i\geq 2$. Clearly,
\begin{eqnarray}\label{proof th1 eq1a}
\Pr\left(\sum_{i=1}^n C_i X_{i,n}>t\right)=\Pr\left(\sum_{i=1}^n C_i X_{i,n}>t,C_1 X_{1,n}>0\right).
\end{eqnarray}
In addition,
\begin{eqnarray}\label{proof th1 eq1c}
\lim_{t\to \infty}\frac{\displaystyle\Pr\left(\sum_{i=1}^n C_i X_{i,n}>t,C_1 X_{1,n}>0\right)}{\Pr(C_1 X_{1,n}\Me{>t})}=\lim_{t\to \infty}\frac{\Pr\bigg(\big(C_1 X_{1,n}/t,C_2 X_{2,n}/t\ldots,C_n X_{n,n}/t\big)\in A_1\bigg)}{\Pr(C_1 X_{1,n}>t)}=\mu(A_1)=1,
\end{eqnarray}
where $A_1:=\{\underline{x}:\sum_{i=1}^n x_i>1, x_1>0, x_i>M \;\mbox{for all}\; i\geq 2\}$ \Me{and $M$ is a negative constant}. Note that the second step is due to the fact that Proposition~A2.12 of Embrechts   et al.\ (1997, p. 563) can be applied in (\ref{Th1-step2- vague conv order stats}) since $A_1$ does not put any mass on its boundary. In addition, the whole mass over the set $A_1$ is concentrated on the line $(1,\infty]\times\{0\}^{n-1}$.
Combining (\ref{proof th1 eq1a})\ignore{, (\ref{proof th1 eq1b})} and (\ref{proof th1 eq1c}) we have
\begin{eqnarray}\label{proof th1 eq1}
\Pr\left(\sum_{i=1}^n C_i X_{i,n}>t\right)\sim\Pr\left(C_1 X_{1,n}>t\right), \quad t\to\infty.
\end{eqnarray}
Similarly,
\begin{eqnarray*}\label{proof th1 eq2a}
\lim_{t\to \infty}\frac{\displaystyle\Pr\left(C_1 X_{1,n}-\sum_{i=2}^n C_i X_{i,n}>t,C_1 X_{1,n}>0\right)}{\Pr(C_1 X_{1,n}\Me{>t})}
&=&\mu(A_2)=1,
\end{eqnarray*}
where $A_2:=\{\underline{x}:x_1-\sum_{i=2}^n x_i>1, x_1>0, x_i>M \;\mbox{for all}\; i\geq 2\}$. Once again, the entire mass over the set $A_2$ is concentrated on the line $(1,\infty]\times\{0\}^{n-1}$. \ignore{Note that
\begin{eqnarray}\label{proof th1 eq2b}
\Pr\left(C_1 X_{1,n}-\sum_{i=2}^n C_i X_{i,n}>t,C_1 X_{1,n}\leq 0\right)=0,
\end{eqnarray}
since $C_i\geq 0$ for all $i\geq 2$.} Thus, \ignore{relations (\ref{proof th1 eq2a}) and (\ref{proof th1 eq2b}) suggest that}
\begin{eqnarray}\label{proof th1 eq2}
\Pr\left(C_1 X_{1,n}-\sum_{i=2}^n C_i X_{i,n}>t\right)\sim \Pr(C_1 X_{1,n}>t), \quad t\to\infty.
\end{eqnarray}
We may now drop the non-negativity assumption for the $C_i, i \geq 2$ since
\BQN \label{TT}
\Pr\left(C_1 X_{1,n}-\sum_{i=2}^n C^{-}_i X_{i,n}>t\right)\leq \Pr\left(\sum_{i=1}^n C_i X_{i,n}>t\right)\leq \Pr\left(\sum_{i=1}^n C^{+}_i X_{i,n}>t\right),
\EQN
where $C^{+}_i=\max\{C_i,0\}$ and $C^{-}_i=\max\{-C_i,0\}$. The latter, together with \hE{(\ref{Th1-step1}),} (\ref{proof th1 eq1}) and  (\ref{proof th1 eq2}) completes the proof for this case. \QED

%
\prooftheo{th3}
\dE{
We show first that for any index $i$ such that $2\le i\le n$
\[ \frac{\hE{C_i^+} X_{i,n}}{t} \Bigl \lvert ( C_1 X_{1,n}> t) \toprob 0
\]
\aE{is valid where} $\toprob$ \jE{and $\todis$ stand for convergence in probability and in distribution, respectively} letting the argument $t \to \infty$.
Indeed, \dE{by \eqref{suff1 joint tail cond}} 
for any $y>0$ we obtain \IME{applying Breiman's Lemma}
\dE{(set $\CMF:=\max_{2 \le i \le n} \hE{C_i^+} $)}
\BQNY
\Pr\left( \frac{\dE{C_i^+ X_{i,n}}}{t}> y \bigl \lvert C_1 X_{1,n}>t\right)
&=& \frac{ \Pr\big(\dE{C_i^+} X_{i,n}> ty, C_1 X_{1,n}>t\big)}{\pk{C_1 X_{1,n}> t}}\\
&\le & \frac{ \Pr\big(\CMF X_{2,n}> ty \big)}{\pk{C_1 X_{1,n}> t}}\\
&\le & \frac{\Pr(X_1>t y)}{\Pr(X_1>t)} \frac{\Pr(X_1>t)}{\Pr(C_1 X_{1,n}>t)}
\sum_{1\leq i<j\leq n}
\frac{ \Pr\big( \CMF X_i>t y , \CMF X_j>t y\big) }{\pk{X_1>t y}}\nonumber\\
&\to & 0\quad \text{ as } t\to \infty
\EQNY
\jE{since by the assumptions on $\CMF$
$$\lim_{t\to 0}\frac{ \Pr\big( \CMF X_i>t y , \CMF X_j>t y\big) }{\pk{X_1>t y}} =0.$$
}
Therefore, $\lim_{t\to \IF} \Pr\left( \frac{\dE{C_i^+ X_{i,n}}}{t}> y \bigl \lvert C_1 X_{1,n}>t\right)=0$. Further, equation (\ref{Th1-step1}) implies
that $\frac{C_1 X_{1,n}}{t} \Bigl \lvert (C_1 X_{1,n}> t) \todis W,$ where \aE{the random variable} $W\ge 1$ \jE{has survival function} $x^{-\alpha}, x\ge 1$.
 Thus,
\BQN\label{ime}
\left(\frac{C_1 X_{1,n}}{t}, \frac{C_2^+ X_{2,n}}{t} , \ldot \frac{\hE{C_n^+} X_{n,n}}{t} \right) \Bigl \lvert ( C_1 X_{1,n}> t)  \todis (W,0, \ldots,0)
\EQN
implying
\BQN
\Biggl( \frac{\IME{C_1} X_{1,n}}{t} + \displaystyle\sum_{i=2}^n \frac{C_i^+  \IME{X_{i,n}}}{t} \Biggr)\Bigl \lvert (C_1 X_{1,n}> t) \todis W.
\EQN
When $\alpha=0$, then $W=1$ and hence the convergence holds in probability. \IME{Similarly, we obtain
$$\Biggl(\frac{\IME{C_1} X_{1,n}}{t} + \displaystyle\sum_{i=2}^n \frac{C_i^- \IME{X_{i,n}}}{t}\Biggr) \Bigl \lvert (C_1 X_{1,n}> t) \todis W.$$
Consequently, }
\def\ve{\varepsilon}
\eE{
\BQNY
 \lim_{t\to \infty} \frac{ \Pr\left( C_1 X_{1,n}- \displaystyle\sum_{i=2}^n C_i^{-} X_{i,n}> t\right)}{\pk{C_1 X_{1,n}> t }} &=&
 \lim_{t\to \infty} \frac{ \Pr\left( C_1 X_{1,n}-  \displaystyle\sum_{i=2}^n C_i^{-} X_{i,n}> t, C_1 X_{1,n}> t \right)}{\pk{C_1 X_{1,n}> t }} \\&=&
 \lim_{t\to \infty} \Pr\left(C_1 X_{1,n}-  \sum_{i=2}^n C_i^{-} X_{i,n}> t \bigl \lvert C_1 X_{1,n}> t  \right) = \pk{W >1}=1.
 \EQNY
If $\CMF = 0$ almost surely the proof follows, therefore let us assume that $\CMF>0$. \hE{Suppose for notational simplicity that $C_i>0, i\le n$.} For any $\ve >0$ we have
\BQNY
\Pr\left(\sum_{i=1}^n C_i X_{i,n}>t\right)
&\le &
\Pr\left(\sum_{i=1}^n C_i X_{i,n}>t, C_1 X_{1,n}>t(1- \ve)\right)+\Pr\left(\eE{\sum_{i=2}^n C_i X_{i,n}}> \ve t\right)\nonumber\\
&\hE{\le} &
\Pr\left(\sum_{i=1}^n C_i X_{i,n}>t ,C_1 X_{1,n}>t(1- \ve) \right)+
\Pr\left(\eE{n \CMF  X_{2,n}}> \ve t\right).
\EQNY
By \eqref{Th1-step2-part i} we have $\lim_{t\to \IF} \frac{ \Pr\big(\eE{n \CMF  X_{2,n}> \ve t\big)}}{\pk{ C_1 X_{1,n}> t}}=0.$
Thus, \IME{in view of \eqref{ime}}
\BQNY
\lim_{t\to \IF} \frac{ \displaystyle\Pr\left(\sum_{i=1}^n C_i X_{i,n}>t,C_1 X_{1,n}>t(1- \ve) \right)}{ \pk{ C_1 X_{1,n}> t(1-\ve)}}
=\lim_{t\to \IF} \Pr\left(\sum_{i=1}^n C_i X_{i,n}>t/(1-\ve) \lvert C_1 X_{1,n}>t\right)= \Pr\big( W>1/(1-\ve)\big),
\EQNY
and hence the proof follows from equation \eqref{TT} and
}
letting $\ve \downarrow 0$. \QED
}


\prooftheo{Th2} The proof is based on \hE{\nelem{Lm2} and} similar arguments \hE{as} provided \hE{for} \netheo{Th1}. We first assume that \hE{$C_i\geq 0$ for all $2\leq i\leq n$ and let $M<0$ such that $-M>(n-1)L^*$, where $L^*=\max_{1\leq i<j\leq n} L_{ij}$. Obviously,
\begin{eqnarray}\label{proof th3 eq1a}
\Pr\left(\sum_{i=1}^n C_i X_{i,n}>t\right)=\Pr\left(\sum_{i=1}^n C_i X_{i,n}>t,C_1 X_{1,n}\leq t+M a(t)\right)+\Pr\left(\sum_{i=1}^n C_i X_{i,n}>t,C_1 X_{1,n}>t+M a(t)\right).
\end{eqnarray}}
\hE{\qE{Further, we have}
\begin{eqnarray}\label{proof th3 eq1b}
\Pr\left(\sum_{i=1}^n C_i X_{i,n}>t,C_1 X_{1,n}\leq t+M a(t)\right)&\leq& \Pr\left(\sum_{i=2}^n C_i X_{i,n}>-M a(t)\right)\\
&\leq&\Pr\left( C^* X_{2,n}>-\frac{M}{n-1} a(t)\right)\nonumber\\
&\leq& \sum_{1\leq i<j\leq n} \Pr\left( C^* X_i>-\frac{M}{n-1} a(t),C^* X_j>-\frac{M}{n-1} a(t)\right)\nonumber\\
&=&o\big(\Pr(C_1 X_1>t)\big), \quad t\to\infty,\nonumber
\end{eqnarray}
where the last implication is due to (\ref{suff2-cond2}) and the fact that $-M>(n-1)L_{ij}$ for all $1\leq i<j\leq n$. Now, the second term from $(\ref{proof th3 eq1a})$ is investigated  via (\ref{Lemma2- vague conv order stats}), i.e.,
\begin{eqnarray}\label{proof th3 eq1c}
\lefteqn{\lim_{t\to \infty}\frac{\displaystyle\Pr\left(\sum_{i=1}^n C_i X_{i,n}>t,C_1 X_{1,n}>t+M a(t)\right)}{\Pr(C_1 X_{1,n}>t)}}\\
&=&\frac{\Pr\Big(\big( (C_1X_{1,n}- t)/a(t),C_2X_{2,n}/a(t)\ldot C_n X_{n,n}/a(t)\big)\in B_1\Big)}{\Pr(C_1 X_{1,n}>t)}=\nu(B_1)=1,
\nonumber
\end{eqnarray}
where} $B_1:=\{\underline{x}:\sum_{i=1}^n x_i>0, x_1>M\}$. Now, the second step is due to the fact that Proposition~A2.12 of Embrechts  et al.\ (1997, p. 563) can be applied in (\ref{Lemma2- vague conv order stats}) since $B_1$ does not put any mass on its boundary. In addition, the whole mass over the set $B_1$ is concentrated on the line $(0,\infty]\times\{0\}^{n-1}$. \hE{Combining equations (\ref{proof th3 eq1a}), (\ref{proof th3 eq1b}) and (\ref{proof th3 eq1c}), we get
\begin{eqnarray}\label{proof th3 eq1d}
\Pr\left(\sum_{i=1}^n C_i X_{i,n}>t\right)\sim\Pr(C_1 X_{1,n}>t), \quad t\to\infty.
\end{eqnarray}}
\hE{Similarly,
\begin{eqnarray}\label{proof th3 eq2a}
\lefteqn{\lim_{t\to \infty}\frac{\displaystyle\Pr\left(C_1 X_{1,n}-\sum_{i=2}^n C_i X_{i,n}>t,C_1 X_{1,n}> t+M a(t)\right)}{\Pr(C_1 X_{1,n}>t)}}\\
&=&\frac{\Pr\Big(\big( (C_1X_{1,n}- t)/a(t),C_2X_{2,n}/a(t)\ldot C_n X_{n,n}/a(t)\big)\in B_2\Big)}{\Pr(C_1 X_{1,n}>t)}=\nu(B_2)=1,
\nonumber
\end{eqnarray}
where} $B_2:=\{\underline{x}:x_1-\sum_{i=2}^n x_i>0, x_i>M\}$. Once again, the entire mass over the set $B_2$ is concentrated on the line $(0,\infty]\times\{0\}^{n-1}$. Note that
\[\Pr\left(C_1 X_{1,n}-\displaystyle\sum_{i=2}^n C_i X_{i,n}>t,C_1 X_{1,n}\leq t+M a(t)\right)=0\]
due to the non-negativity assumption of the $C_i$'s. The latter and (\ref{proof th3 eq2a}) yield that
\begin{eqnarray}\label{proof th3 eq2b}
\Pr\left(C_1 X_{1,n}-\sum_{i=2}^n C_i X_{i,n}>t\right)\sim\Pr(C_1 X_{1,n}>t), \quad t\to\infty.
\end{eqnarray}
Therefore, equations (\ref{proof th3 eq1d}), (\ref{proof th3 eq2b}) and (\ref{TT}) help in dropping the non-negativity assumption for the $C_i$ for all $2\leq i \leq n$, and the proof follows
utilising further \eqref{LS}.\QED

\pE{\textbf{Acknowledgments}: We are in debt to two reviewers for numerous comments which improved the manuscript.
The second author kindly acknowledges partial support by an SNF grant.}

\end{document}